\documentclass{amsart}
\usepackage{amssymb, amsmath,  amsfonts}
\usepackage{color}
\usepackage{slashbox, longtable}
\usepackage[all,cmtip]{xy}

\newtheorem{thm}{Theorem}[section]

\newtheorem{example}[thm]{Example}

\newcommand{\circline}{\!\;\!\!\circ\!\!\;\!\!-\!\!\!-\!}

\numberwithin{equation}{section}

\input{tcilatex}

\begin{document}
\title[An update of $QH^*(G/P)$]{An update of quantum cohomology of
homogeneous varieties}
\author{Naichung Conan Leung}
\address{The Institute of Mathematical Sciences and Department of
Mathematics, The Chinese University of Hong Kong, Shatin, Hong Kong}
\email{leung@math.cuhk.edu.hk}
\author{Changzheng Li}
\address{Kavli Institute for the Physics and Mathematics of the Universe
(WPI), The University of Tokyo, 5-1-5 Kashiwa-no-Ha,Kashiwa City, Chiba
277-8583, Japan}
\email{changzheng.li@ipmu.jp}

\begin{abstract}
We describe recent progress on $QH^{\ast }\left( G/P\right) $ with special
emphasis of our own work.
\end{abstract}

\maketitle

{\allowdisplaybreaks[4] }





\section{Introduction}


How many intersection points are there for two given lines in a plane? How
many lines pass through two given points in a plane? It is the main concern
in enumerative geometry to find solutions to such questions of counting
numbers of geometric objects that satisfy certain geometric conditions.
There are two issues here. First, we should impose conditions so that the
expected solution to a counting problem is a finite number. Second, we will
work in the complex projective space, so that Schubert's \textit{principle
of conservation of number} holds. Then we ask for this invariant. For
instance to either of the aforementioned questions, the solution is a
constant number $1$, if the condition is imposed precisely as ``two distinct
given lines (resp. points) in a complex projective plane $\mathbb{P}^2$".

Information on counting numbers of geometric objects may be packaged to form
an algebra. For instance for the case of a complex projective line $\mathbb{P%
}^1$, basic geometric objects are either a point or the line $\mathbb{P}^1$
itself. There are only two pieces of non-trivial enumerative information
among three basic geometric objects. Namely, (a) $\langle \mbox{pt}, %
\mbox{line}, \mbox{line}\rangle_{0, 3, 0}=1$, telling us that a point and
two (same) lines intersect at a unique intersection point; (b) $\langle %
\mbox{pt}, \mbox{pt}, \mbox{pt}\rangle_{0, 3, 1}=1$, telling us that there
is a unique line passing through three (distinct) given points.
Incorporating all these enumerative information together, we obtain the
algebra $\mathbb{C}[x, q]/\langle x^2-q\rangle$. Here the identity $1$ and
the element $x$ stand for the basic geometric objects, a line and a point,
respectively; $q$ stands for a line, the geometric object to be counted; the
aforementioned counting numbers are then read off directly from the
algebraic relations $1* x=x$ and $x* x= q$, in terms of the coefficient $%
N_{1, x}^{x, 0}$ of $x$ and the coefficient $N_{x, x}^{1, q}$ of $q$,
respectively. In modern language, we are saying that the quantum cohomology
ring $QH^*(\mathbb{P}^1)$ of $\mathbb{P}^1$ is isomorphic to $\mathbb{C}[x,
q]/\langle x^2-q\rangle$ as an algebra.

The concept of quantum cohomology of a smooth projective complex variety
arose from the subject of string theory in theoretic physics in 1990s, and
the terminology was introduced by the physicists \cite{Vafa}. The
coefficients for the quantum multiplication are genus zero Gromov-Witten
invariants, which were rigorously defined later (by means of virtual
fundamental classes in general) via symplectic geometry \cite{RuTi} and via
algebraic geometry \cite{KoMa}. As a first surprising application of the
\textit{big} quantum cohomology, Kontsevich solved an old problem in
enumerative geometry on counting the number of rational curves of degree $d$
passing through $3d-1$ points in general position in $\mathbb{P}^{2}$, by
giving a recursive formula in 1994. We note that a rational curve of degree $%
1$ is a complex projective line. The space $\mathbb{P}^{1}$ can be written
as the quotient of the special linear group $SL(2,\mathbb{C})$ by its
subgroup of upper triangular matrices $\Big(%
\begin{array}{cc}
a & b \\
0 & a^{-1} \\
&
\end{array}%
\Big)$, and it is a special case of the so-called homogeneous variety $G/P$.
In the present article, we will be concerned only with the (\textit{small})
quantum cohomology ring $QH^{\ast }(G/P)$ of 
$G/P$, which is an algebra deforming the classical cohomology ring $H^{\ast
}(G/P)$ by incorporating all genus zero, three-point \textit{Gromov-Witten
invariants} as the coefficients for the (small) quantum multiplication.
Namely, it consists of all information on counting numbers of rational
curves of various (fixed) degrees that pass through three given projective
subvarieties. The study of the ring structure of $QH^{\ast }(G/P)$ is
referred to as the quantum Schubert calculus, which is not only a branch of
algebraic geometry, but also of great interest in mathematical physics,
algebraic combinatorics, representation theory, and so on.

In the present article, we will give a brief review of the developments on
quantum Schubert calculus since Fulton's beautiful lecture \cite%
{Fulton-survey}, with a focus on the authors' work \cite%
{LeungLi-functorialproperties, LeungLi-GWinv, LeungLi-QuantumToClassical,
LeungLi-QPR, Li-functorialproperties}. We apologize to the many whose work
less related with the four problems listed in section \ref{subset-QSC} has
not been cited, and apologize to those whose work we did not notice.

\section{A brief review of Schubert calculus and generalizations}

\subsection{Classical Schubert calculus for homogeneous varieties $G/P$}

Like Grimm's fairy tales for the children, the next enumerative problem is
known to everybody in the world of Schubert calculus:

\textit{How many lines in $\mathbb{P}^3$ intersect four given lines in
general position?}

\noindent The solution to it is $2$, in modern language, obtained by a
calculation in the cohomology ring of the complex Grassmannian $Gr(2, 4)$
(see e.g. \cite{KlLa}).

A complex Grassmannian $Gr(k,n)$ consists of $k$-dimensional vector
subspaces in $\mathbb{C}^{n}$, i.e., one-step flags $V\leqslant \mathbb{C}%
^{n}$. In particular, 
$\mathbb{P}^{n-1}=Gr(1,n)$. A direct generalization of it is the variety of
partial flags
\begin{equation*}
F\ell _{n_{1},\cdots ,n_{r};n}:=\{V_{n_{1}}\leqslant \cdots \leqslant
V_{n_{r}}\leqslant \mathbb{C}^{n}~|~\dim V_{n_{i}}=n_{i},\forall 1\leq i\leq
r\},
\end{equation*}%
where $[n_{1},\cdots ,n_{r}]$ is a fixed subsequence of $[1,2,\cdots ,n-1]$.
Every partial flag variety $F\ell _{n_{1},\cdots ,n_{r};n}$ is a quotient $%
SL(n,\mathbb{C})/P$ of $SL(n,\mathbb{C})$ by one parabolic subgroup $P$ that
consists of block-upper triangular matrices of the following type:
\begin{equation*}
\left(
\begin{array}{cccc}
M_{1} & \ast  & \ast  & \ast  \\
\mathbf{0} & M_{2} & \ast  & \ast  \\
\vdots  & \ddots  & \ddots  & \vdots  \\
\mathbf{0} & \cdots  & \mathbf{0} & M_{r+1}%
\end{array}%
\right) _{n\times n},
\end{equation*}%
where $M_{i}$ is an $(n_{i}-n_{i-1})\times (n_{i}-n_{i-1})$ invertible matrix  for $1\leq i\leq
r$, and $n_{0}:=0,n_{r+1}:=n$. These are called homogeneous varieties (or
flag varieties) of Lie type $A_{n-1}$.

In general, homogeneous varieties are $X=G/P$, where $G$ is a
simply-connected complex simple Lie group, and $P$ a parabolic subgroup.
These are classified by data from Dynkin diagrams, i.e., by a pair $(\Delta
,\Delta _{P})$ of
sets $\Delta _{P}\subset \Delta $. In particular, the parabolic subgroup
corresponding to the empty subset of $\Delta $ is called a Borel subgroup,
denoted as $B$ instead. Let $W$ denote the Weyl group of $G$, and $W_{P}$
denote the Weyl subgroup corresponding to $P$. There is always an
accompanied combinatorial subset $W^{P}\subset W$ bijective to $W/W_{P}$,
which consists of minimal length representatives of the cosets in $W/W_{P}$
with respect to a canonical length function $\ell :W\rightarrow \mathbb{Z}%
_{\geq 0}$. Every $\sigma ^{u}$ is 
the class of a Schubert variety $\Omega _{u}$ of complex codimension $\ell
(u)$, and therefore is in $H^{2\ell (u)}(X,\mathbb{Z})$. The coefficients $%
N_{u,v}^{w}$ in the cup product in $H^{\ast }(X)$,
\begin{equation*}
\sigma ^{u}\cup \sigma ^{v}=\sum_{w}N_{u,v}^{w}\sigma ^{w},
\end{equation*}%
count the number of intersection points of three Schubert varieties $g\Omega
_{u},g^{\prime }\Omega _{v},g^{\prime \prime }\Omega _{w^{\sharp }}$ with
generic elements $g,g^{\prime },g^{\prime \prime }\in G$, where $w^{\sharp
}\in W^{P}$ parameterizes the dual basis of $H^{\ast }(G/P)$ to $\{\sigma
^{u}\}$ with respect to the bilinear form $(\sigma ^{u},\sigma
^{v}):=\int_{[G/P]}\sigma ^{u}\cup \sigma ^{v}$. In particular, the
structure constants $N_{u,v}^{w}$ are nonnegative integers, which are known
as Littlewood-Richardson coefficients in the special case of complex
Grassmannians. Those unfamiliar with the general theory may just note the
following two points, and then refer to \cite{Fulton-survey} for a very nice
introduction to $Gr(k,n)$.

(1) When $G/P=F\ell_{n_1, \cdots, n_r; n}$, the Weyl group $W$ of $G$ is the
group $S_n$ of permutations of $\{1, \cdots, n\}$. The subset $\Delta_P$ is $%
\Delta\setminus\{\alpha_{n_1}, \cdots, \alpha_{n_r}\}$, provided that $%
\Delta=\{\alpha_1, \cdots, \alpha_{n-1}\}$ is a base of simple roots of Lie$%
(G)$ whose Dynkin diagram is canonical, i.e., given by
\begin{tabular}{l}
\raisebox{-0.4ex}[0pt]{$  \circline\!\;\!\!\circ\cdots\,
\circline\!\!\!\;\circ $} \\
\raisebox{1.1ex}[0pt]{${\hspace{-0.2cm}\scriptstyle{\alpha_1}\hspace{0.3cm}%
\alpha_2\hspace{0.8cm}\alpha_{n-1} } $}%
\end{tabular}%
\!. In particular, $G/B=F\ell_{1, 2, \cdots, n-1; n}$.

(2) When $G/P=Gr(k, n)$, the combinatorial subset $W^P$ is given by
\begin{equation*}
W^{P}:=\{w\in S_n~|~ w(1)<w(2)<\cdots< w(k); w(k+1)<w(k+2)<\cdots< w(n)\}.
\end{equation*}
Take a partition $\lambda=(\lambda_1, \cdots, \lambda_k)$ with $n-k\geq
\lambda_1\geq \lambda_2\geq \cdots\geq \lambda_k\geq 0$, and a complete flag
$F_\bullet$ of vector subspaces $\{0\}=:F_0\leqslant F_1 \leqslant\cdots \leqslant
F_{n-1}\leqslant F_n:=\mathbb{C}^n$ with $\dim_{\mathbb{C}}F_i=i$ for $0\leq
i\leq n$. We have the Schubert variety
\begin{equation*}
\Omega_\lambda(F_\bullet):=\{V\leqslant \mathbb{C}^n~|~ \dim (V\cap
F_{n-k+i-\lambda_i})\geq i, \forall 1\leq i\leq k\}\subset Gr(k, n).
\end{equation*}
Then $N_{u, v}^w$ counts the number of intersection points of three Schubert
varieties $\Omega_{\lambda(u)}(F_\bullet)$, $\Omega_{\lambda(v)}(F_\bullet^{%
\prime})$, $\Omega_{\lambda(w^\sharp)}(F_\bullet^{\prime\prime})$ with
respect to three general flags and the partitions $\lambda(u)=(u(k)-k,
\cdots, u(2)-2, u(1)-1)$, $\lambda(v)=(v(k)-k, \cdots, v(2)-2, v(1)-1)$ and $%
\lambda(w^\sharp)=(n-k+1-w(1), n-k+2-w(2), \cdots, n-w(k))$.

A thorough study of 
$H^{\ast }(X)$ includes (but not limited to) a combinatorial description of
the Littlewood-Richardson coefficients, a ring presentation of $H^{\ast }(X)$
with certain generators, and an expression of every Schubert class $\sigma
^{u}$ in terms of polynomial in the aforementioned generators. In the case
of complex Grassmannians, there have been nice answers to all the above
through various approaches. Much is also known about $H^{\ast }(G/P)$. However, a \textit{manifestly positive} combinatorial rule for the
structure constants $N_{u,v}^{w}$ is still lacking in general. To the
authors' knowledge, such rules have been shown only for two-step flag
varieties $F\ell _{n_{1},n_{2};n}$ \cite{Coskun-twostep} (and there is a
preprint \cite{Coskun-partialflag}), besides (co)minuscule Grassmannians.
One may have noted that every Schubert variety admits a resolution by an
associated Bott-Samelson manifold, which is a tower of $\mathbb{P}^{1}$%
-fibrations. Using the topology of the Bott-Samelson resolution, Duan
obtained a nice algorithm for computing $N_{u,v}^{w}$ for general $G/P$ in
\cite{Duan}, though it involves sign cancellation.

\subsection{Equivariant Schubert calculus}

Every parabolic subgroup $P$ of $G$ contains a Borel subgroup $B$. Let $K$
be a maximal compact Lie subgroup of $G$, such that $T=K\cap B$ is a maximal
torus of $K$. We consider the $T$-equivariant cohomology $H_{T}^{\ast
}(G/P)=H_{T}^{\ast }(G/P,\mathbb{Q})$, which is a module over the ring $%
S:=H_{T}^{\ast }(\mbox{pt})$ which is a polynomial ring $\mathbb{Q}[\alpha
_{1},\cdots ,\alpha _{n}]$ freely generated by simple roots $\alpha _{j}$'s
of $G$. $H_{T}^{\ast }(G/P)$ has an $S$-basis of Schubert classes $\sigma
^{u}$'s, indexed by the same combinatorial set $W^{P}$. The structure
constants $p_{u,v}^{w}$'s in the equivariant product in $H_{T}^{\ast }(G/P)$%
,
\begin{equation*}
\sigma ^{u}\circ \sigma ^{v}=\sum_{w}p_{u,v}^{w}\sigma ^{w},
\end{equation*}%
are homogeneous polynomials of degree $\ell (u)+\ell (v)-\ell (w)$ in $S$
(by which we mean the zero polynomial if the degree is negative). The
evaluation of the polynomial $p_{u,v}^{w}(\alpha _{1},\cdots ,\alpha _{n})$
at the origin coincides with the intersection number $N_{u,v}^{w}$.
Geometric meanings of $p_{u,v}^{w}$'s are not completely clear.
Nevertheless, $p_{u,v}^{w}$ enjoys the positivity property, lying inside $%
\mathbb{Q}_{\geq 0}[\alpha _{1},\cdots ,\alpha _{n}]$ (or enjoys a form of positivity conjectured by Peterson and proved by Graham \cite{Grah},
  as an element in $\mathbb{Q}_{\geq 0}[-\alpha _{1},\cdots
,-\alpha _{n}]$, up to the choice of Schubert varieties that determine the
equivariant Schubert classes ). The study of $H_{T}^{\ast }(G/P)$
is referred to as the \emph{equivariant Schubert calculus}.

\subsection{Affine Schubert calculus}

The \emph{affine Kac-Moody group} $\mathcal{G}$ associated to $G$ is an
algebraic analog of the free (smooth) loop space $LK=Map\left( \mathbb{S}%
^{1},K~\right) $. The Dynkin diagram of $\mathcal{G}$ is precisely the
affine Dynkin diagram of $G$, namely it has an extra simple root, called the
\emph{affine root} $\alpha _{0}$, attached to the simple root corresponding
to the adjoint representation of $G$. In particular $\Delta _{\mathrm{aff}%
}=\Delta \cup \{\alpha _{0}\}$ is a base of $\mathcal{G}$.

A natural generalization of $G/P$ is the infinite dimensional projective
ind-variety $\mathcal{G}/\mathcal{P}$, where $\mathcal{P}$ denotes a
parabolic subgroup of $\mathcal{G}$. The two extreme cases $\mathcal{B}$, $%
\mathcal{P}_{{\scriptsize \mbox{max}}}$ correspond to the subsets $\emptyset
$, $\Delta $ of $\Delta _{\mathrm{aff}}$ respectively\footnote{%
If the subset of $\Delta _{\mathrm{aff}}$ contains $\alpha _{0}$, then $%
\mathcal{G}/\mathcal{P}$ reduces to $G/P$.
\par
{}}. The corresponding ind-varieties $\mathcal{G}/\mathcal{B}$ and $\mathcal{%
G}/\mathcal{P}_{{\scriptsize \mbox{max}}}$ are called \emph{affine flag
manifolds} and \emph{affine Grassmannians }respectively. Indeed $\mathcal{G}/%
\mathcal{B}$ is homotopy equivalent to $LK/T$ and $\mathcal{G}/\mathcal{P}_{%
{\scriptsize \mbox{max}}}$ is homotopy equivalent to the \emph{base loop
group }$\Omega K=\{f\in LK~|~f(\mbox{id}_{\mathbb{S}^{1}})=\mbox{id}_{K}\}$.
Both of them inherit the $T$-action. One may study the $T$-equivariant
cohomology rings of $LK/T$ and $\Omega K$. Both of them admit an $S$-basis
of affine Schubert classes.

On the other hand the (equivariant) homology $H_{\ast }^{T}(\Omega K)$ (in
the sense of Borel-Moore homology) also admits a natural ring structure,
called the \emph{Pontryagin product }which is induced by the group
multiplication of $K$. The space $H_{\ast }^{T}(\Omega K)$ is further
equipped with a Hopf algebra structure, with the coproduct structure induced
from the (equivariant) cohomology ring structure of $\Omega K$. The study of
this structure is referred to as the (equivariant) \emph{affine Schubert
calculus}. A surprising fact, as will be discussed in the next section, is
that $H_{\ast }^{T}(\Omega K)$ is essentially the same as the equivariant
\emph{quantum cohomology} of $G/B$.\qquad

\subsection{Quantum Schubert calculus}

\label{subset-QSC} The (small) \emph{quantum cohomology }ring $QH^{\ast
}(G/P)$ is part of the intersection theory on the moduli space of stable
maps to $G/P$ (see e.g. \cite{FuPa}). It contains the ordinary $H^{\ast
}(G/P)$ as $G/P$ is a connected component of the moduli space of stable maps
to $G/P$, namely the component of constant maps. $QH^{\ast }(G/P)=(H^{\ast
}(G/P)\otimes \mathbb{Q}[q_{1},\cdots ,q_{r}],\ast )$ also has a basis of
Schubert classes $\sigma ^{u}$ over $\mathbb{Q}[q_{1},\cdots ,q_{r}]$, where
$r=\dim H_{2}(G/P)$. We identify $H_{2}(G/P,\mathbb{Z})$ with $\mathbb{Z}%
^{r} $ with basis given by two (real) dimensional Schubert cycles in $G/P$.

Given $\mathbf{d}=(d_{1},\cdots ,d_{r})\in \mathbb{Z}^{r}=H_{2}(G/P,\mathbb{Z%
})$, we denote $q^{\mathbf{d}}:=q_{1}^{d_{1}}\cdots q_{r}^{d_{r}}$. The
structure constants $N_{u,v}^{w,\mathbf{d}}$ in the quantum multiplication,
\begin{equation*}
\sigma ^{u}\ast \sigma ^{v}=\sum_{w\in W^{P},\mathbf{d}\in H_{2}(G/P,\mathbb{%
Z})}N_{u,v}^{w,\mathbf{d}}\sigma ^{w}q^{\mathbf{d}},
\end{equation*}%
are genus zero, three-point Gromov-Witten invariants. Geometrically, $%
N_{u,v}^{w,\mathbf{d}}$ counts the cardinality of the set
\begin{equation*}
\{f:\mathbb{P}^{1}\rightarrow G/P~|~f_{\ast }([\mathbb{P}^{1}])=\mathbf{d}%
;f(0)\in g\Omega _{u},f(1)\in g^{\prime }\Omega _{v},f(\infty )\in g^{\prime
\prime }\Omega _{w^{\sharp }};f\mbox{ is holomorphic}\}
\end{equation*}%
with generic $g,g^{\prime },g^{\prime \prime }\in G$. In particular, $%
N_{u,v}^{w,\mathbf{d}}\in \mathbb{Z}_{\geq 0}$, and it is zero unless $%
d_{i}\geq 0$ for all $i$. A holomorphic map $f:\mathbb{P}^{1}\rightarrow G/P$
of degree $\mathbf{0}$ is a constant map. Therefore $N_{u,v}^{w,\mathbf{0}%
}=N_{u,v}^{w}$ for the cup product.

The study of $QH^{\ast }(G/P)$ is referred to as the \emph{quantum Schubert
calculus}, which includes at least the following  as pointed out in \cite%
{Fulton-survey}.

\begin{enumerate}
\item A presentation of the ring, $QH^{\ast }(G/P)=\mathbb{Q}[x_{1},\cdots
,x_{N},q_{1},\cdots ,q_{r}]/(\mathrm{relations})$. See section 3.1.

\item A \textit{manifestly positive} combinatorial rule for the structure
constants $N_{u,v}^{w,\mathbf{d}}$.

\item A \textquotedblleft quantum Giambelli" formula, which expresses each
Schubert class $\sigma ^{u}$ as a polynomial in the generators $x_{i}$ and $%
q_{j}$. See section 3.4.

\quad Because of the lack of functoriality for \emph{quantum} products in
general, we also have the following problem.

\item A comparison between $QH^{\ast }(G/B)$ and $QH^{\ast }(G/P)$. See
section 3.2.
\end{enumerate}

We can consider the $T$-equivariant version $QH_{T}^{\ast }(G/P)$ of the
quantum cohomology in a similar fashion, which is a module over $%
S[q_{1},\cdots ,q_{r}]$ generated by Schubert classes $\sigma ^{u}$'s. The
structure constants $\tilde{N}_{u,v}^{w,\mathbf{d}}$'s in the equivariant
quantum product
\begin{equation*}
\sigma ^{u}\star \sigma ^{v}=\sum_{w,\mathbf{d}}\tilde{N}_{u,v}^{w,\mathbf{d}%
}\sigma ^{w}q^{\mathbf{d}},
\end{equation*}%
are again homogeneous polynomials in $S$, enjoying the positivity property
(provided that a choice of positive/negative simple roots is chosen
properly) \cite{Miha-positivity}. They contain information on both $QH^{\ast
}(G/P)$ and $H_{T}^{\ast }(G/P)$: (i) $\tilde{N}_{u,v}^{w,\mathbf{0}}$
coincides with $p_{u,v}^{w}$ for the $T$-equivariant product $\sigma
^{u}\circ \sigma ^{v}$ and (ii) the evaluation of the polynomial $\tilde{N}%
_{u,v}^{w,\mathbf{d}}\in S$ at the origin coincides with $N_{u,v}^{w,\mathbf{%
d}}$ for the quantum product.

We remark that it is possible to define $QH_{T}^{\ast }(G/P)$ over $\mathbb{Z%
}$, instead of $\mathbb{C}$.

\section{An overview on (equivariant) quantum Schubert calculus}

\subsection{Ring presentations}

\label{subset-ringpresentation} In order to find a ring presentation of $%
QH^{\ast }(G/P)$, it is natural to start with an explicit ring presentation
of the ordinary cohomology ring $H^{\ast }\left( G/P\right) $, say%
\begin{equation*}
H^{\ast }(G/P)\cong {\frac{\mathbb{Q}[x_{1},\cdots ,x_{N}]}{\langle f_{1}(%
\mathbf{x}),\cdots ,f_{m}(\mathbf{x})\rangle }.}
\end{equation*}%
Such a presentation is know in many cases. Then a lemma of Siebert and Tian
\cite{SiTi} (or Proposition 11 of \cite{FuPa}) shows that the ring structure
on $QH^{\ast }(G/P)=H^{\ast }(G/P)\otimes \mathbb{Q}[q_{1},\cdots ,q_{r}]$
can be obtained by deforming the relations $f_{i}$'s, i.e.%
\begin{equation*}
QH^{\ast }(G/P)\cong {\frac{\mathbb{Q}[x_{1},\cdots ,x_{N},q_{1},\cdots
,q_{r}]}{\langle f_{1}(\mathbf{x})-g_{1}(\mathbf{x},\mathbf{q}),\cdots
,f_{m}(\mathbf{x})-g_{m}(\mathbf{x},\mathbf{q})\rangle }.}
\end{equation*}%
Here each $g_{i}(\mathbf{x},\mathbf{q})$ is computed by $f_{i}$ after
replacing the cup product by the quantum product. This is usually difficult
to compute.

For instance, $H^{\ast }(\mathbb{P}^{1})\overset{\cong }{\rightarrow }%
\mathbb{Q}[x]/\langle x^{2}\rangle $, which sends the hyperplane class $H$
(Poincar\'{e} dual to $[\mbox{pt}]$) to $x$. Since $H\ast H=q$ in $QH^{\ast
}(\mathbb{P}^{1})$, we have $QH^{\ast }(\mathbb{P}^{1})\cong \mathbb{Q}%
[x]/\langle x^{2}-q\rangle $.

Most known ring presentations of $QH^{\ast }(G/P)$ are obtained in this way.
Such $G/P$'s include $F\ell _{n_{1},\cdots ,n_{r};n}$ \cite{GiKi},\cite{
Kim-partialflag},\cite{AsSa} and various Grassmannians \cite{SiTi},\cite{
KrTa-Lagr},\cite{KrTa-Orth},\cite{CMP-RingPres},\cite{BKT-Isotropic},\cite%
{ChPe-Adjoint}. Some Grassmannians of exceptional Lie types  are still unknown.

\bigskip

There is another way to get a ring presentation of $QH^{\ast }(G/P)$, by
finding \textquotedblleft quantum differential equations". Givental's $J$%
-function is a $H^{\ast }(G/P)$-valued function, involving gravitational
correlators (a class of invariants more general than Gromov-Witten
invariants). It was introduced for any smooth projective variety $X$, and
played an important role in mirror symmetry. Quantum differential equations
are certain differential operators annihilating $J$. Every quantum
differential equation gives rise to a relation in $QH^{\ast }(X)$. (See e.g.
Example 10.3.1.1 of \cite{CoKa-book} for the case of $\mathbb{P}^{1}$). With
this method, Kim studied the quantum $D$-module of $G/B$ and obtained the
following ring presentation of $QH^{\ast }(G/B)$ for general $G$.

\begin{thm}[\protect\cite{Kim-GoverB}]
The small quantum cohomology ring $QH^{\ast }(G/B,\mathbb{C})$ is
canonically isomorphic to
\begin{equation*}
\mathbb{C}[p_{1},\cdots ,p_{l},q_{1},\cdots ,q_{l}]/I,
\end{equation*}%
where $l$ denotes the rank of $G$, and $I$ is the ideal generated by the
nonconstant complete integrals of motions of the Toda lattice for the
Langlands-dual Lie group $(G^{\vee },B^{\vee },(T^{\mathbb{C}})^{\vee })$ of
$(G,B,T^{\mathbb{C}})$.
\end{thm}

\noindent   There have been the descriptions of $J$-function for general $G/P$   \cite%
{BCK-GW-AbNonAb},\cite{Braverman}. Nevertheless, even for complex
Grassmannians, there are no closed formulas on the quantum differential
equations, to the authors' knowledge.

\bigskip

In his unpublished lecture notes \cite{Peterson}, Dale Peterson announced a
uniform presentation of $QH^{\ast }(G/P)$ (and its $T$-equivariant
extension) for all $P$. The so-called \emph{Peterson variety} $Y$ is
subvariety in $G^{\vee }/B^{\vee }$, which is equipped with a $\mathbb{C}%
^{\ast }$-action. The $\mathbb{C}^{\ast }$-fixed points $y_{P}$'s are
isolated and parameterized by the finite set of (conjugacy classes of)
parabolic subgroups $P$. (See \cite{Rietsch-petersoniso},\cite{HaTy} for
precise descriptions for type $A$ case.) By considering $y\in Y$ with $%
z\cdot y$ approaching $y_{P}$ with various $P$'s as $z\in \mathbb{C}^{\ast }$
goes to $0$ or infinity, we obtain two stratifications of $Y$ by affine
varieties $Y_{P}^{+}$'s or $Y_{P}^{-}$'s respectively.

Peterson claimed that the spectra of the quantum cohomology ring $QH^{\ast
}(G/P)$ is $Y_{P}^{+}$, or equivalently $QH^{\ast }(G/P)\cong \mathbb{C}%
[Y_{P}^{+}]$ as algebras. This was proved by Rietsch \cite%
{Rietsch-petersoniso} for all $SL(n,\mathbb{C})/P$, and by Cheong \cite%
{Cheong} when $G/P$ is a Lagrangian Grassmannian or an orthogonal
Grassmannian.

\bigskip

Peterson also interpreted all $\mathbb{C}[Y_{P}^{-}]$ as the homology of
based loop groups. In particular, $Y_{G}^{-}$ is birational to $Y_{B}^{+}$
and $\mathbb{C}[Y_{G}^{-}]\cong H_{\ast }(\Omega K)$, where $K$ is a maximal
compact subgroup of $G$. Its consequence has become the following theorem
now, as was firstly proved by Lam and Shimozono. There is also an
alternative proof by the authors.

\begin{thm}[\protect\cite{Peterson},\protect\cite{LaSh-GoverPaffineGr},%
\protect\cite{LeungLi-GWinv}]
\label{QHandHomology} The equivariant quantum cohomology ring $QH_{T}^{\ast
}(G/B)$ is isomorphic to the equivariant homology $H_{\ast }^{T}(\Omega K)$
as algebras, after localization.
\end{thm}

\noindent The above isomorphism is explicit, sending a (localized) $S$-basis
$\sigma ^{w}q^{\mathbf{d}}$ of $QH_{T}^{\ast }(G/B)$ to a (localized) $S$%
-basis of Schubert homology classes of $H_{\ast }^{T}(\Omega K)$. The proof in \cite{LaSh-GoverPaffineGr} also showed
that $QH_{T}^{\ast }(G/P)$ is isomorphic to a quotient of $H_{\ast
}^{T}(\Omega K)$, after localization.

The algebraic structures of $QH_{T}^{\ast }(G/B)$, or more generally $%
QH_{T}^{\ast }(G/P)$, is completely determined by the equivariant quantum
Chevalley formula (i.e. quantum multiplication by elements in $H^{2}(G/P)$
in section \ref{subset-combrule}) together with a few natural properties.
This criterion was obtained by Mihalcea \cite{Miha-EQCRandCri}. Using this
criterion, the proof of the above theorem can be reduced to explicit
computations of certain products in $H_{\ast }^{T}(\Omega K)$. These
products were first obtained by Lam and Shimozono, mainly by using
Peterson's $j$-isomorphism \cite{Lam} together with properties of Kostant
and Kumar's nilHecke algebras \cite{KoKu}. Later we realized that they can
also be obtained by carefully analyzing a combinatorial formula on the
structure constants (to be described in section \ref{subset-combrule}).
Despite the two aforementioned (combinatorial) proofs, a satisfactory
understanding of the above theorem is still lacking.

When $G=SL(n,\mathbb{C})$, the connections among Givental-Kim's presentation
\cite{GiKi,Kim-GoverB}, Peterson's presentation above, and Kostant's
solution to Toda lattice \cite{Kostant}, are now becoming better understood
\cite{LaSh-QspTodalattice},\cite{LaSh-DoubleQspTodalattice}.
 In this case, the quantum cohomology
  $QH^{\ast }(SL(n,\mathbb{C})/B)$ is   also closely related with Dunkl elements, which   leads to some relevant applications \cite{Kiri}.
Some other characterizations for $QH^{\ast }(G/B)$ can be found for example in  \cite{Mare-relations},%
\cite{ Mare-characterization}.

Finally, we remark that  there are some studies on the quantum differential equations for the cotangent bundle of a complete flag variety $G/B$ \cite{BMO} or of a    partial flag variety  $F\ell_{n_1, \cdots, n_r; n}$  \cite{MaOk}, \cite{GRTV}, \cite{TaVa}. This might lead to nice applications to the quantum cohomology of the corresponding flag variety by taking an appropriate limit. For instance,    an application  to the Chevalley-type formula was given  in \cite{BMO} (as will be discussed in   section \ref{sub-combrule}). The    cotangent bundle of a partial flag variety is a Nakajima quiver variety.
The quantum cohomology of it can be related with the integral systems and quantum groups \cite{NeSh},\cite{BMO},\cite{MaOk}.
 The relation between such an approach and the aforementioned Peterson's approach for complex Grassmannians is studied in   \cite{GoKo}.

 \subsection{Comparing $QH^{\ast }(G/P)$ with $QH^{\ast }(G/B)$}

\label{subset-comparison copy(1)}The inclusions $B\subset P\subset G$ induce
a fiber bundle
\begin{equation*}
\pi :G/B\rightarrow G/P
\end{equation*}
with fiber $P/B$, which is again a complete flag variety.

\begin{example}
\label{example-FL3} For $G=SL(3,\mathbb{C})$ with $B\subsetneq P\subsetneq G$%
, we have $G/B=F\ell _{1,2;3}$
and $G/P=Gr(2,3)\cong \mathbb{P}^{2}$; $\pi $ coincides with the natural
forgetful map, sending a flag $V_{1}\leqslant V_{2}\leqslant \mathbb{C}^{3}$
in $G/B$ to the partial flag $V_{2}\leqslant \mathbb{C}^{3}$ in $G/P$; the
fiber of $\pi $ is $P/B=\mathbb{P}^{1}$.
\end{example}

For ordinary cohomologies, the Leray-Serre spectral sequence relates the
cohomology ring of $G/B$ with those of $G/P$ and $P/B$. For instance the
induced homomorphism $\pi ^{\ast }:H^{\ast }(G/P)\rightarrow H^{\ast }(G/B)$
is injective, sending the Schubert class $\sigma _{P}^{u}$ for $G/P$ (where $%
u$ varies over the combinatorial subset $W^{P}\subset W=W^{B}$) to the
Schubert class $\sigma _{B}^{u}$ for $G/B$.

For quantum cohomologies, there is no such functoriality in general. For
instance, the pullback map $\pi ^{\ast }$ on $H^{\ast }$ does not even have
a \emph{quantum }analog. Nevertheless, such functoriality does exist in our
specific case. The analog of $\pi ^{\ast }$ is a comparison formula stated
by Peterson \cite{Peterson} and proved by Woodward \cite{Wood}.

\begin{thm}[Peterson-Woodward comparison formula]
\label{comparison} Every structure constant $N_{u,v}^{w, \lambda_P }$ for
the quantum product $\sigma^u_P*\sigma^v_P$ in $QH^*(G/P)$ coincides with a
structure constant $N_{u,v}^{w\tilde{w}, \lambda_B}$ for
$\sigma^u_B*\sigma^v_B$ in $QH^*(G/B)$, where $(\tilde w, \lambda_B)\in
W_P\times H_2(G/B, \mathbb{Z})$ is uniquely and explicitly determined by $%
\lambda_P\in H_2(G/P, \mathbb{Z})$.
\end{thm}

\noindent We remark that the equivariant quantum extension of the above
comparison formula is implicitly contained in Corollary 10.22 of \cite%
{LaSh-GoverPaffineGr} (see \cite{LiHu-EQSC} for more details).

For ordinary cohomologies, $H^{\ast }(G/B)$ admits a $\mathbb{Z}^{2}$%
-filtration $\mathcal{F}$, which induces an isomorphism $Gr^{\mathcal{F}%
}(H^{\ast }(G/B))\cong H^{\ast }(P/B)\otimes H^{\ast }(G/P)$ of $\mathbb{Z}%
^{2}$-graded algebras, by the spectral sequence. In order to generalize this
to quantum cohomologies, we need to find a nice grading on $QH^{\ast }(G/B)$%
. This is quite tricky indeed.

\noindent \textbf{Example \ref{example-FL3}. (continued) }{\itshape The Weyl
group $W$ of $SL(3,\mathbb{C})$ is the permutation group $S_{3}$, generated
by transpositions $s_{1}=(12)$ and $s_{2}=(23)$. The cohomology degree of a
Schubert class $\sigma ^{u}$ is $2\ell (u)$.\ Explicitly $\ell (u)=$}$%
0,1,1,2,2,3$ when $u=1,s_{1},s_{2},s_{1}s_{2},s_{2}s_{1},s_{1}s_{2}s_{1}$
respectively. {The $\mathbb{Z}^{2}$-grading $gr(\sigma ^{u})$ of $\sigma
^{u} $}$\in ${$H^{\ast }(F\ell _{1,2;3})$, as given by the spectral sequence
for }$\pi :${$F\ell _{1,2;3}\rightarrow Gr(2,3)$, are $%
(0,0),(1,0),(0,1),(0,2),(1,1),(1,2)$ respectively. This grading gives a
number of nice consequences on $H^{\ast }(F\ell _{1,2;3})$. }

The quantum cohomology $QH^{\ast }(F\ell _{1,2;3})=H^{\ast }(F\ell
_{1,2;3})\otimes \mathbb{Q}[q_{1},q_{2}]$ has a $\mathbb{Q}$-basis $\sigma
^{w}q_{1}^{a}q_{2}^{b}$. %
The above $\mathbb{Z}^{2}$-grading map extends to $QH^{\ast }(F\ell
_{1,2;3}) $, by defining
\begin{equation*}
gr(\sigma ^{u}q_{1}^{a}q_{2}^{b}):=gr(\sigma ^{u})+(2a-b,3b).
\end{equation*}%
It is tricky to find this grading, e.g. $gr(q_{2})=(-1,3)$. This is the
\emph{correct }grading as all consequences from spectral sequence have
natural quantum generalizations for $QH^{\ast }(F\ell _{1,2;3})$. (See
Example 1.1 of \cite{LeungLi-functorialproperties} for precise descriptions.)

\vspace{0.15cm}

In general, we need to take a \textit{maximal} chain of parabolic subgroups $%
P_{j}$, i.e.,
\begin{equation*}
B:=P_{0}\subsetneq P_{1}\subsetneq \cdots \subsetneq P_{r-1}\subsetneq
P_{r}=P\subsetneq G,
\end{equation*}
where $r$ is the rank of the Levi subgroup of $P$. This corresponds to a
chain of subsets $\emptyset :=\Delta _{0}\subsetneq \Delta _{1}\subsetneq
\cdots \subsetneq \Delta _{r-1}\subsetneq \Delta _{r}=\Delta _{P}$ with $%
|\Delta _{i}|=i$.

We can always find $P_{j}$'s such that $P_{j}/P_{j-1}$'s are all projective
spaces $\mathbb{P}^{N_{j}}$, with at most one exception occurring at the
last step $P_{r}/P_{r-1}$. For instance, for $P\subset Sp(8,\mathbb{C})$
with $Sp(8,\mathbb{C})/P\simeq \mathbb{P}^{7}$, we have $r=3$ with $P_{1}/B=%
\mathbb{P}^{1}$, $P_{2}/P_{1}=\mathbb{P}^{2}$ and $P/P_{2}=LGr(3,6)$ is a
Lagrangian Grassmannian. Precise choices are made in Table 2 of \cite%
{LeungLi-functorialproperties} when $P/B$ is of type $A$ (and its associated Dynkin diagram is
connected), and in section 3.5 of \cite{LeungLi-functorialproperties} (or
Table 1 of \cite{Li-functorialproperties}) for a general $P/B$.

The key point of the whole story, is to find a nice $\mathbb{Z}^{r+1}$%
-grading on $QH^{\ast }(G/B)$ with respect to the chosen chain. The
Peterson-Woodward comparison formula plays a key role in defining such a
grading\footnote{%
The original definition of this grading was made recursively in \cite%
{LeungLi-functorialproperties}. It was greatly simplified in \cite%
{Li-functorialproperties}, by proving a conjecture due to a referee of \cite%
{LeungLi-functorialproperties}.}. With this grading, the authors obtained
certain functorial properties among quantum cohomologies of homogeneous
varieties in the following sense\footnote{%
We remark that part of the statements were only proved for $P/B$ of type $A$
in \cite{LeungLi-functorialproperties}, and are proved for all general cases
in \cite{Li-functorialproperties} recently.}.

\begin{thm}[\protect\cite{LeungLi-functorialproperties}, \protect\cite%
{Li-functorialproperties}]
\label{functorial} Let $\pi: G/B\to G/P$ denote the natural projection, and $%
r$ denote the rank of the Levi subgroup of $P$.

\begin{enumerate}
\item There exists a $\mathbb{Z}^{r+1}$-filtration $\mathcal{F}$ on $%
QH^*(G/B)$,
respecting the quantum product structure. 

\item There exist an ideal $\mathcal{I}$ 
of $QH^*(G/B)$ and a canonical algebra isomorphism
\begin{equation*}
QH^*(G/B)/\mathcal{I} \overset{\simeq}{\longrightarrow}QH^*(P/B).
\end{equation*}

\item There exists a subalgebra $\mathcal{A}$ of $QH^*(G/B)$ together with
an ideal $\mathcal{J}$ of $\mathcal{A}$,
such that $QH^*(G/P)$ is canonically isomorphic to $\mathcal{A}/\mathcal{J}$
as algebras.

\item There exists a canonical injective morphism of graded algebras
\begin{equation*}
\Psi_{r+1}:\,\, \,\, QH^*(G/P)\hookrightarrow Gr^{\mathcal{F}%
}_{(r+1)}\subset Gr^{\mathcal{F}}(QH^*(G/B))
\end{equation*}
together with an isomorphism of graded algebras after localization
\begin{equation*}
Gr^{\mathcal{F}}(QH^*(G/B))\cong\big(\bigotimes_{j=1}^{r}QH^*(P_{j}/P_{j-1})%
\big)\bigotimes Gr^{\mathcal{F}}_{(r+1)},
\end{equation*}
where $P_j$'s are parabolic subgroups constructed in a canonical way,
forming a chain $B:=P_0\subsetneq P_1\subsetneq\cdots\subsetneq
P_{r-1}\subsetneq P_r=P\subsetneq G$.

Furthermore, $\Psi_{r+1}$ is an isomorphism if and only if $P_{j}/P_{j-1}$
is a projective space for any $1\leq j\leq r$.
\end{enumerate}

All the relevant ideals, subalgebras and morphisms above could be described
concretely.
\end{thm}

When only quantum multiplication with an element in $H^{2}(G/B)$ is
involved, this theorem is reduced to the quantum Chevalley formula. To prove
the theorem, we use induction with respect to the length $\ell (u)$ of
Schubert classes $\sigma ^{u}$. The Peterson-Woodward comparison formula is
used frequently and the positivity of the structure constants turns out to
be needed as well. The most dedicated arguments occur in the proof of (1)
and in the part to show $\Psi _{r+1}$ is an morphism for the general case.
The notion of virtual coroot was introduced to reduce many cases in general
Lie types to type $A$, while there are still a number of cases that require
individual discussions.

The functorial properties of $QH^{\ast }(G/P)$, as given in theorems \ref%
{comparison} and \ref{functorial}, have many applications in finding
combinatorial rules on certain $N_{u,v}^{w,\mathbf{d}}$'s, especially on the
so-called \emph{quantum to classical} \emph{principle}. We could also relate
certain Gromov-Witten invariants between $G/P$ and $G^{\prime }/P^{\prime }$
for $G\neq G^{\prime }$. We will discuss these next.

\subsection{Combinatorial rules}

\label{subset-combrule} The problem of finding a manifestly positive
combinatorial rule for the Gromov-Witten invariants, or equivalently the
structure constants, $N_{u,v}^{w,\mathbf{d}}$ of $QH^{\ast }(G/P)$ is open.
Even the special case when $\mathbf{d}=0$, namely the counterpart for
classical cohomology $H^{\ast }(G/P)$, is still not solved except in very
limited cases, including complex Grassmannians.

\subsubsection{A combinatorial formula on $N_{u,v}^{w,\mathbf{d}}$ with
signs cancellations}

If we do not require coefficients of the combinatorial rule to be all
positive, namely we allow signs cancellations, then it suffices to find one
for $QH^{\ast }(G/B)$ because of the Peterson-Woodward comparison formula.
Using the relationship between $QH^{\ast }\left( G/B\right) $ and $H_{\ast
}(\Omega K)$, we obtained such a formula as we explain next.

Let $\Delta =\{\alpha _{1},\cdots ,\alpha _{n}\}\subset \mathfrak{h}^{\ast }$
denote a base of simple roots of $G$ and $\rho $ (resp. $\rho ^{\vee }$)
denote the sum of fundamental (co)weights.
Then $H_{2}(G/B,\mathbb{Z})$ is canonically identified with the coroot
lattice $Q^{\vee }:=\bigoplus_{i=1}^{n}\mathbb{Z}\alpha _{i}^{\vee }$. For
the affine flag manifold $LK/T$, $H_{T}^{\ast }(LK/T)$ has an $S$-basis of
Schubert classes $\mathfrak{S}^{x}$, with $x=wt_{\lambda }$, parametrized by
the affine Weyl group $W_{{\scriptsize \mbox{af}}}=W\ltimes Q^{\vee }$. Here
$t_{\lambda }$ denotes the translation by $\lambda $.

To describe the combinatorial rule, we define the Kostant (homogeneous)
polynomial $d_{y,[x]}$ and a rational function $c_{x,[y]}$ as below: (i) $%
d_{x,y}$ is the coefficient $p_{x,y}^{y}$ of $\mathfrak{S}^{y}$ in $%
\mathfrak{S}^{x}\cdot \mathfrak{S}^{y}=\sum p_{x,y}^{z}\mathfrak{S}^{z}$;
(ii) $c_{x,[y]}=\sum_{w\in W}c_{x,yw}$ with $\big(c_{x^{\prime },y^{\prime }}%
\big)$ being the transpose inverse of the (infinitely dimensional
invertible) matrix $\big(d_{x,y}\big)$. They both have combinatorial
descriptions \cite{KoKu, Kumar}.

\begin{thm}[\protect\cite{LeungLi-GWinv}]
\label{formulawithSign}The structure constants $N_{u,v}^{w,\lambda }$ in the
quantum product of $QH^{\ast }\left( G/B\right) $,
\begin{equation*}
\sigma ^{u}\ast \sigma ^{v}=\sum_{w\in W,\lambda \in Q^{\vee
}}N_{u,v}^{w,\lambda }\sigma ^{w}q^{\lambda },
\end{equation*}%
%
are given by the constant function
\begin{equation*}
N_{u,v}^{w,\lambda }=\sum_{\lambda _{1},\lambda _{2}\in Q^{\vee
}}c_{ut_{A},[t_{\lambda _{1}}]}c_{vt_{A},[t_{\lambda
_{2}}]}d_{wt_{2A+\lambda },[t_{\lambda _{1}+\lambda _{2}}]}
\end{equation*}%
where $A:=-12n(n+1)\rho ^{\vee }$, provided that $\langle 2\rho ,\lambda
\rangle =\ell (u)+\ell (v)-\ell (w)$ and $\lambda \succcurlyeq 0$, and zero
otherwise.
\end{thm}

\noindent The proof of this theorem uses the fact that $H_{\ast }^{T}(\Omega
K)$ has a basis $\psi _{t_{\lambda }}$ parametrized by translation
representatives of $W_{{\scriptsize \mbox{af}}}/W$ for which the Pontryagin
product is simple: $\psi _{t_{\lambda _{1}}}\cdot \psi _{t_{\lambda
_{2}}}=\psi _{t_{\lambda _{1}+\lambda _{2}}}$; the change of bases to
Schubert classes is known explicitly by localization and the relationship
between $QH^{\ast }\left( G/B\right) $ and $H_{\ast }(\Omega K)$.

\subsubsection{Combinatorial formulae of Pieri-Chevalley type}\label{sub-combrule}

There are manifestly positive combinatorial formulas for quantum
multiplication by \emph{special} Schubert classes $\sigma ^{v}$. When $\ell
(v)=1$, i.e. $v=s_{i}$ is a simple refection in $W$, it is called the
Chevalley formula for $H^{\ast }(G/P)$. The quantum version of Chevalley
formula for $QH^{\ast }(G/P)$ was conjectured by Peterson \cite{Peterson},
and was first proved by Fulton and Woodward \cite{FuWo}. It is also called
the quantum Monk formula in the special case of $G/P=F\ell _{1,\cdots ,n-1;n}
$ \cite{FGP}, and it has an alternative description when $G/P$ is a
(co)minuscule Grassmannian \cite{CMP-RingPres}.

The $T$-equivariant generalization was also conjectured by Peterson \cite%
{Peterson}. It was first proved by Mihalcea \cite{Miha-equivquantumGrass,
Miha-EQCRandCri} in a combinatorial way, and there is also a geometric
approach \cite{BuMi-curveNbd} by the technique of curve neighborhoods. A
special case of it is stated as follows, for which there is another proof by
using the Springer resolution \cite{BMO}.

\begin{thm}[Equivariant quantum Chevalley formula for $G/B$]
\label{quantumChevalley}In $QH_T^*(G/B)$,
\begin{equation*}
\sigma^u\star \sigma^{s_i}=(\chi_i-u(\chi_i))\sigma^u+ \sum \langle \chi_i,
\gamma^\vee \rangle \sigma^{us_\gamma}+ \sum \langle \chi_i, \gamma^\vee
\rangle q^{\gamma^\vee}\sigma^{us_\gamma},
\end{equation*}
where $\chi_i$ denotes a fundamental weight, the first sum is over positive
roots $\gamma$ with $\ell(us_\gamma)=\ell(u)+1$, and the second sum is over
positive roots $\gamma$ with $\ell(us_\gamma)=\ell(u)+1-\langle 2\rho,
\gamma^\vee\rangle$.
\end{thm}

\noindent The corresponding formula for $G/P$ is slightly more complicated.

\bigskip

When $G/P$ is a complex Grassmannian $Gr(k,n)$, there is an exact sequence
of tautological bundles over it: $0\rightarrow \mathcal{S}\rightarrow
\underline{\mathbb{C}}^{n}\rightarrow \mathcal{Q}\rightarrow 0$. The fiber
of $\mathcal{S}$ over $[V]\in Gr(k,n)$ is given by the vector subspace $V$
in $\mathbb{C}^{n}$ itself. The ordinary cohomology $H^{\ast }(Gr(k,n))$, as
a ring, is generated by Chern classes $c_{p}(\mathcal{S})$'s (or $c_{p}(%
\mathcal{Q})$'s). A (manifestly positive) combinatorial formula on the
multiplication by either sets of Chern classes is referred to as a \emph{%
Pieri rule}. The quantum version of a Pieri rule was first given by Bertram
\cite{Bertram}.

The analog of $Gr\left( k,n\right) $ for other Lie types is $G/P_{%
{\scriptsize \mbox{max}}}$ with $P_{{\scriptsize \mbox{max}}}$ a maximal
parabolic subgroup of $G$. When $G$ is a classical group, $G/P_{{\scriptsize %
\mbox{max}}}$ parameterizes linear subspaces which are isotropic with
respect to a non-degenerate bilinear form which is skew-symmetric (for type $%
C$) or symmetric (for type $B$ and $D$). Therefore such a $G/P_{{\scriptsize %
\mbox{max}}}$ is usually called an \emph{isotropic Grassmannian}. The
corresponding quantum Pieri rules with respect to Chern classes of
tautological quotient bundles have been obtained by Buch, Kresch and
Tamvakis \cite{KrTa-Lagr, KrTa-Orth}, \cite{BKT-Isotropic}.

For instance when $G/P=IG(k,2n)$ is a (non-maximal) isotropic Grassmannian
of type $C$, Schubert classes can also labeled by $(n-k)$-strict partitions,
with $c_{p}(\mathcal{Q})$ corresponding to a class $\sigma ^{p}$ of the
special $(n-k)$-strict partition $p$. In terms of $(n-k)$-strict partitions,
one has a quantum Pieri rule with respect to $c_{p}(\mathcal{Q})$'s in the
following form.

\begin{thm}[\protect\cite{BKT-Isotropic}]
\begin{equation*}
\sigma^p * \sigma^\lambda=\sum_\mu2^{N(\lambda,
\mu)}\sigma^\mu+\sum_\nu2^{N(\lambda, \nu^\sharp)-1}\sigma^{\nu} q,
\end{equation*}
where $\nu^{\sharp}$ is an $(n-k)$-strict partition for $IG(k+1, 2n+2)$,
associated to the $(n-k)$-partition $\nu$ for $IG(k, 2n)$.
\end{thm}

\noindent The classical part of the above formula is new even for the
classical cohomology $H^{\ast }(IG(k,2n))$.

The quantum Pieri rules with respect to $c_{p}(\mathcal{S})$'s have been
studied by the authors in \cite{LeungLi-QPR}. For instance, when $%
G/P=IG(k,2n)$, there is another parameterization of the Schubert classes by
\textit{shapes}, which are pairs of partitions. In terms of shapes, every
Chern class $c_{p}(\mathcal{S}^{\ast })$ corresponds to a class $\sigma ^{p}$
of special shape $p$, and one has a quantum Pieri rule with respect to these
Chern classes in the following form.

\begin{thm}[\protect\cite{LeungLi-QPR}]
\begin{equation*}
\sigma^{p}\star\sigma^{\mathbf{a}}=\sum_{\mathbf{b}}2^{e(\mathbf{a},\mathbf{b%
})}\sigma^{\mathbf{b}}+\sum_{\mathbf{c}}2^{e(\tilde{\mathbf{a}},\tilde{%
\mathbf{c}})}\sigma^{\mathbf{c}}q,
\end{equation*}
where $\tilde{\mathbf{a}}$ and $\tilde{\mathbf{c}}$ are shapes for $IG(k-1,
2n)$, 
associated to the shapes $\mathbf{a}$ and $\mathbf{c}$ for $IG(k, 2n)$
respectively.
\end{thm}

\noindent The classical part of the above formula is the classical Pieri
rule of Pragacz and Ratajski \cite{PrRa-Pierirule}. We remark that when $G/P$
is a non-maximal isotropic Grassmanian of type $B$ or $D$, the above formula
does involve sign cancellations even for some degree \emph{one}
Gromov-Witten invariants, thus it is not quite satisfactory.

\bigskip

For homogeneous varieties of type $A$, namely partial flag varieties, there
are natural forgetting maps to complex Grassmannians $\pi _{n_{i}}:F\ell
_{n_{1},\cdots ,n_{r};n}\rightarrow Gr(n_{i},n)$. The ring $QH^{\ast }(F\ell
_{n_{1},\cdots ,n_{r};n})$ is generated by the pull-back of Chern classes of
the tautological subbundle (or quotient bundle) over $Gr(n_{i},n)$ for all $%
i $. The quantum Pieri rule with respect to these classes was obtained by
Ciocan-Fontanine \cite{CFon}. The equivariant quantum version of it has been
obtained in \cite{LiHu-EQSC} recently.

All these quantum Pieri rules are obtained by determining relevant
Gromov-Witten invariants $N_{u,v}^{w,\mathbf{d}}$ of $QH^{\ast }(G/P)$
explicitly.

\subsubsection{Calculations of $N_{u, v}^{w, \mathbf{d}}$}

Each $N_{u,v}^{w,\mathbf{d}}$ is an intersection number of cycles in the
moduli space of stable maps, which is the stable maps compactification of
the space of morphisms from $\mathbb{P}^{1}$ to $G/P$. Bertram \cite{Bertram}
used a different compactification, namely the quot schemes compactification,
in order to calculate the relevant Gromov-Witten invariants. This method was
further used in \cite{CFon}, \cite{KrTa-Lagr, KrTa-Orth}.

For $Gr(k,n)$, because of $H_{2}=\mathbb{Z}$, $N_{u,v}^{w,d}$ counts the
number of rational curves $C$ of degree $d$ passing through three Schubert
subvarieties of $Gr(k,n)$. Buch \cite{Buch} introduced the \textit{span}
(resp. \textit{kernel}) of $C$ as the smallest (resp. largest) subspace of $%
\mathbb{C}^{n}$ containing (resp. contained in) all the $k$-dimensional
subspaces parametrized by points of $C$. Buch showed that the span (resp.
kernel) of $C$ determines a point in a related Schubert subvariety of $%
Gr(k+d,n)$ (resp. $Gr(k-d,n)$). Furthermore, he related Gromov-Witten
invariants involved in a quantum Pieri rule for $Gr(k,n)$ to classical
intersection numbers (of Pieri type) of $Gr(k+1,n)$, thus giving an
elementary proof \cite{Buch}.

This idea was later used by Buch, Kresch and Tamvakis to show that all
Gromov-Witten invariants $N_{u,v}^{w,d}$ for complex Grassmannians,
Lagrangian Grassmannians and (maximal) orthogonal Grassmannians are classical. Such a
phenomenon is now referred to as the \emph{quantum to classical} principle.
It was further shown to hold for the remaining two (co)minuscule
Grassmannians of exceptional Lie type, i.e., Cayley plane $E_{6}/P_{1}$ and
the Freudenthal variety $E_{7}/P_{7}$. Combining both statements, we have

\begin{thm}[\protect\cite{BKT-quantumtoclassical}, \protect\cite%
{CMP-RingPres}]
All genus zero, three-point Gromov-Witten invariants for a (co)minuscule
Grassmannian $G/P$ are equal to classical intersection numbers on some
auxiliary homogeneous varieties of $G$ .
\end{thm}
\noindent Such a statement has been extended to the equivariant quantum K-theory setting in \cite{BuMi-quantumKtheory}, \cite{ChPe-Rationality}.
The above theorem  leads to a manifestly positive combinatorial
formula for all $N_{u,v}^{w,d}$ for $Gr(k, n)$, because of the known
positive formula on the classical intersection numbers on two-step partial
flag varieties \cite{Coskun-twostep}.  The kernel-span technique was also
used in \cite{BKT-Isotropic} to derive the quantum Pieri rules with respect
to $c_{p}(\mathcal{Q})$'s for non-maximal isotropic Grassmannians.

There is another (combinatorial) approach to show the \textquotedblleft
quantum to classical" principle by the authors \cite%
{LeungLi-QuantumToClassical}: As a consequence of the special case of
Theorem \ref{functorial} with $\pi :G/B\rightarrow G/P$ being a $\mathbb{P}%
^{1}$-bundle, the authors obtained vanishing and identities among various
Gromov-Witten invariants. For any simple root $\alpha $, we introduce a map $%
\mbox{sgn}_{\alpha }:W\rightarrow \{0,1\}$ defined by $\mbox{sgn}_{\alpha
}(w):=1$ if $\ell (w)-\ell (ws_{\alpha })>0$, and $0$ otherwise.

\begin{thm}[\protect\cite{LeungLi-QuantumToClassical}]
\label{quantumtoclassical} For any $u,v,w\in W$ and for any $\lambda \in
Q^{\vee }\simeq H_{2}(G/B,\mathbb{Z})$, we have the following for $QH^{\ast
}\left( G/B\right) $

\begin{enumerate}
\item $N_{u, v}^{w, \lambda}=0$ unless {\upshape $\mbox{sgn}%
_\alpha(w)+\langle \alpha, \lambda\rangle \leq \mbox{sgn}_\alpha(u)+%
\mbox{sgn}_\alpha(v)$} for all $\alpha\in \Delta.$

\item Suppose {\upshape $\mbox{sgn}_\alpha(w)+\langle \alpha, \lambda\rangle
=\mbox{sgn}_\alpha(u)+\mbox{sgn}_\alpha(v)=2$} for some $\alpha\in \Delta$,
then {\upshape
\begin{equation*}
N_{u, v}^{w, \lambda}=N_{us_\alpha, vs_\alpha}^{w, \lambda-\alpha^\vee}=
\begin{cases}
N_{u, vs_\alpha}^{ws_\alpha, \lambda-\alpha^\vee}, & i\!f \mbox{ sgn}%
_\alpha(w)=0 \\
\vspace{-0.3cm} &  \\
N_{u, vs_\alpha}^{ws_\alpha, \lambda}, & i\!f \mbox{ sgn}_\alpha(w)=1 \,\,
{}_{\displaystyle .}%
\end{cases}%
\end{equation*}
}
\end{enumerate}
\end{thm}

Combining the above theorem with the Peterson-Woodward comparison formula,
the authors obtained a lot of nice applications, including the quantum Pieri
rules with respect to $c_{p}(\mathcal{S})$'s for isotropic Grassmannians
\cite{LeungLi-QPR}. There are $T$-equivariant generalization \cite{LiHu-EQSC}
of these results, giving nice applications on equivariant quantum Schubert
calculus, including equivariant quantum Pieri rules for all $SL(n,\mathbb{C}%
)/P$'s.

There are some other ways to compute (part of) $N_{u,v}^{w,\mathbf{d}}$ in a
few cases (see \cite{AmGu}, \cite{BCK-HoriVafaConj}, \cite%
{Post-affineapproach}, \cite{KoSt}, \cite{KLS}, \cite{GoMa} etc.). For
instance for $Gr(k, n)$, the structure constants $N_{u,v}^{w, {d}}$ can also
be computed from the classical intersection numbers on $Gr(k, 2n)$ \cite%
{BCFF}. The equivariant quantum situation of this method is studied in \cite%
{BBT}. For two-step flag variety $F\ell_{n_1, n_2; n}$, the quantum to
classical principle holds for certain $N_{u,v}^{w,\mathbf{d}}$'s as well
\cite{Coskun-symm}.

\subsection{Quantum Giambelli formulae}

These are formulae to express each Schubert class $\sigma ^{u}$ as a
polynomial in the generators $x_{i}$ and $q_{j}$ in a presentation $QH^{\ast
}(G/P)=\mathbb{Q}[\mathbf{x},\mathbf{q}]/(\mbox{relations})$. There are
limited results as even such a ring presentation is still open in many
cases, as discussed in section \ref{subset-ringpresentation}.

Schubert classes of $H^{\ast }(Gr(k,n))$ are labeled by partitions.
Partitions corresponding $c_{p}(\mathcal{S}^{\ast })$'s (or $c_{p}(\mathcal{Q%
})$'s) are called special. The Giambelli formula expresses Schubert classes
in terms of determinants with special Schubert classes as entries. The first
quantum version of Giambelli formula, due to Bertram \cite{Bertram}, was
obtained by evaluation in the classical cohomology ring of quot schemes.
Similar ideas were applied to the case of Lagrangian Grassmannians and
orthogonal Grassmannians \cite{KrTa-Lagr, KrTa-Orth}. There is an
alternative way to obtain a quantum Giambelli formula, by using a quantum
Pieri rule and the classical Giambelli formula. This was used to reprove a
quantum Giambelli formula for complex Grassmannians by Buch \cite{Buch}, and
to obtain one for non-maximal isotropic Grassmannian by Buch, Kresch and
Tamvakis \cite{BKT-quantumGiambelli}. There are also known formulas for some
Grassmannians of exceptional types by using software \cite%
{ChPe-quantumGiambelli}.

For complete flag variety $SL(n,\mathbb{C})/B$,  Formin, Gelfand, and
Postnikov constructed quantum Schubert polynomials  \cite{FGP} that represent the quantum Schubert classes.   For partial flag varieties $%
SL(n,\mathbb{C})/P$, the quantum Giambelli formula are obtained by Ciocan-Fontanine \cite{CFon} by
using a geometric tool of moving lemma for quot scheme. In general, there could be a third way to get a quantum Giambelli formula, by studying the equivariant quantum version first and then taking the non-equivariant limit. The equivariant quantum cohomology ring $QH^*_T(G/P)$ contains more information than the quantum cohomology, but behaves more simply than $QH^*(G/P)$ in the sense that it is essentially  determined by the equivariant quantum Chevalley formula due to a criterion by Mihalcea  \cite{Miha-EQCRandCri}. Suppose that  we already have an expectation on a ring presentation of $QH^*_T(G/P)$, together with an expected formula  on the equivariant quantum Schubert classes. Then we can prove our expectation by checking that these candidate satisfy the equivariant quantum Chevalley formula. In this way,  Mihalcea obtained the equivariant quantum Giambelli formula for   complex
Grassmannians \cite{Miha-equivquantumGiambelli}, and Ikeda,
Mihalcea and Naruse have achieved the case of   maximal isotropic Grassmannians recently  \cite{IMN-EQSC}.
For $SL(n, \mathbb{C})/B$, there are  quantum double Schubert polynomials studied by Kirillov and Maeno \cite{KiMa}, and by Ciocan-Fontanine and Fulton \cite{CFF}. Lam and Shimozono define the analogues for $SL(n, \mathbb{C})/P$ and show them to represent the equivariant quantum Schubert classes by using the third approach
\cite{LaSh-quantumdouble}. Such a result has also been   independently obtained  by Anderson and Chen by deriving an
equivariant moving lemma for quot schemes \cite{AnCh}.

\subsection{A few remarks}

So far we have mainly focused on an overview of the developments on the four
problems listed in section \ref{subset-QSC}. There are many other
interesting problems in quantum Schubert calculus, for instance, the study
of the symmetry among (part of) $QH^{\ast }(G/P)$ or its $T$-equivariant
extension \cite{Belk}, \cite{Post-symm}, \cite{CMP-hiddensymm},\cite{CMP-semisimple},
\cite{CMP-affinesymm}, \cite{Coskun-symm}. While there are fewer results on $%
QH_{T}^{\ast }(G/P)$, especially on combinatorial rules of the structure
constants and on equivariant quantum Giambelli formulas, our discussions are
not made systematic.

On the other hand, because of Theorem \ref{QHandHomology}, (equivariant)
quantum Schubert calculus is essentially part of (equivariant) affine
Schubert calculus. For instance, (equivariant) quantum Pieri rules \cite%
{LMSS}, \cite{LaSh-affinePieri} can be obtained this way. We apologize for
not mentioning all such applications.

A natural generalization of (equivariant) quantum Schubert calculus is
(equivariant) \emph{quantum K-theory} of homogeneous varieties \cite%
{Givental-quantumK}, \cite{Lee-quantumK}. However, very little is known,
including \cite{BuMi-quantumKtheory}, \cite{ChPe-Rationality}, \cite%
{BCMP-finiteness}, \cite{LiMi}, etc. Peterson's approach to the homology of
affine Grassmannians could be generalized to obtain a K-theoretic analogue
\cite{LSS-KaffineSC}. It is interesting to know:

\vspace{0.15cm}

\textit{Is there a K-theoretic analogue of Theorem \ref{QHandHomology}, or
more generally, of the statements on strata data $Y_{P}^{\pm }$ of the
Peterson variety?\footnote{%
Recently, Lam, Mihalcea, Shimozono and the second author have made such a
conjectural K-theoretic analogue.}}

\vspace{0.15cm}

There are other generalizations of the quantum Schubert calculus, say for
some \textit{inhomogeneous} varieties, for instance odd symplectic
Grassmannians \cite{Pech}.

\bigskip

The notion of quantum cohomology arose in string theory in mathematical
physics. It is natural to ask

\vspace{0.15cm} 
\textit{What is mirror symmetry of quantum Schubert calculus? }
\vspace{0.15cm}

\noindent Mirror symmetry predicts that the quantum cohomology ring $%
QH^{\ast }(G/P)$ is isomorphic to the Jacobian ring $\mbox{Jac}(W)$ of a
mirror Landau-Ginzberg model $(X^{\vee },W)$. Recall that a \emph{%
Landau-Ginzberg model} is pair $(X^{\vee },W)$, consisting a non-compact K%
\"{a}hler manifold $X^{\vee }$ together with a holomorphic function $%
W:X^{\vee }\rightarrow \mathbb{C}$, which is called a \emph{superpotential}.
For instance when $G/P=\mathbb{P}^{1}$, we have $X^{\vee }=\mathbb{C}^{\ast
} $ and $W:\mathbb{C}^{\ast }\rightarrow \mathbb{C}$ is defined by $W(z)=z+{%
\frac{q}{z}}$, where the quantum parameter $q$ is treated as a fixed nonzero
complex number. Then we have $QH^{\ast }(\mathbb{P}^{1},\mathbb{C})=\mathbb{C%
}[x]/\langle x^{2}-q\rangle \cong \mathbb{C}[z,z^{-1}]/\langle 1-{\frac{q}{%
z^{2}}}\rangle =\mbox{Jac}(W)$.

As we have already seen, such a mirror statement is about a ring
presentation of $QH^{\ast }(G/P)$. Using the presentation of $QH^{\ast
}(G/P) $ announced by Peterson, Rietsch constructed a mirror Landau-Ginzberg
model of $G/P$ \cite{Rietsch-mirrorsymm}. We refer our readers to \cite%
{Lam-Whit} for some recent developments in the relations between the quantum
Schubert calculus and mirror symmetry (as well as the Whittaker functions)
from an algebro-combinatorial perspective.

Mirror symmetry also predicts an isomorphism between $QH^{\ast }(G/P)$ and $%
\mbox{Jac}(W)$ as Frobenius manifolds, matching flat coordinates of both
sides. It would be very interesting if one could get combinatorial rules on $%
N_{u,v}^{w,\mathbf{d}}$ by using the corresponding basis of $\mbox{Jac}(W)$.
In certain cases, Schubert classes do play special roles in mirror symmetry
\cite{GoSm}, \cite{PeRi-LagGrass, PeRi-oddQuadr},\cite{MaRi}.

\section*{Acknowledgements}

The authors would like to thank
Anatol N. Kirillov, Thomas Lam, Leonardo C. Mihalcea and Vijay Ravikumar   for useful discussions and valuable comments.
The first author is supported in part by a RGC grant no. 401411. The second author is   supported in part by JSPS Grant-in-Aid for Young Scientists (B) No. 25870175.

\end{document}